\documentclass{amsart}

\newtheorem{thm}{Theorem}

\newtheorem{lem}[thm]{Lemma}

\newtheorem{prop}[thm]{Proposition}

\theoremstyle{definition}

\newtheorem{exmp}[thm]{Example}

\newtheorem{ack}{Acknowledgments}

\newtheorem{defn-thm}[thm]{Definition--Theorem}  

\theoremstyle{remark}

\renewcommand{\c}[0]{{\mathbb C}}

\newcommand{\hl}{\\\hline}

\newcommand{\lcm}[0]{\operatorname{lcm}}

\def\fract#1#2{\raise4pt\hbox{$ #1 \atop #2 $}}

\def\bfa{{\bf a}}

\def\bfw{{\bf w}}

\def\calo{{\mathcal O}}

\def\calo{{\mathcal O}}

\def\cals{{\mathcal S}}

\def\bbc{{\mathbb C}}

\def\bbp{{\mathbb P}}
\def\bbq{{\mathbb Q}}
\def\bbr{{\mathbb R}}

\def\bbz{{\mathbb Z}}

\def\gre{\epsilon}

\def\gsp1{{\mathfrak s}{\mathfrak p}(1)}

\def\la#1{\hbox to #1pc{\leftarrowfill}}
\def\ra#1{\hbox to #1pc{\rightarrowfill}}
\def\Se{Sasakian-Einstein }
\def\Ke{K\"ahler-Einstein }

\begin{document}
\bibliographystyle{amsalpha}

\title{Einstein metrics on rational homology spheres}
\author{Charles P. Boyer and Krzysztof Galicki}
\address{Department of Mathematics and Statistics, 
University of New Mexico,
Albuquerque, NM 87131.}
\email{cboyer@math.unm.edu}
\email{galicki@math.unm.edu}

\maketitle

\section{Introduction} 

In this paper we prove the existence of Einstein metrics, actually \Se metrics, on nontrivial 
rational homology spheres in all odd dimensions greater than $3.$ It appears as though little is 
known about the existence of Einstein metrics on rational homology spheres, and the known 
ones are typically homogeneous. The are two exception known to the authors.
Both involve Sasakian geometry and both occur in dimension 7.  In \cite{BGN02} the 
two authors and M. Nakamaye gave a list of 184 rational homology 
7-spheres with \Se metrics. The result was based on a theorem of Johnson and Koll\'ar proving
the existance of K\"ahler-Einstein metric on cartain Fano 3-folds with orbifold singularities
\cite{jk2}. More recently, Grove, Wilking and Ziller constructed infinitely many
rational homology 7-spheres with 3-Sasakian metrics of cohomogeneity 
one under an action of $S^3\times S^3$ \cite{GWZ03}. Being 3-Sasakian
these metrics are necessarily Sasakian-Einstein.
Their construction involves orbifold self-dual Einstein 
metrics discovered by N. Hitchin \cite{Hi93} a decade ago.
Hitchin constructed a family of positive self-dual Einstein metrics, indexed by integers $k\geq 4$ 
which live on $S^4\setminus\bbr\bbp^2$. They are complete in the orbifold sense, i.e., they can
be viewed as metrics defined on compact Riemannian orbifolds $M^4_k$
with a $\bbz_k$ quotient singularities along $\bbr\bbp^2$. Any positive self-dual Einstein orbifold 
$\calo$ admits a compact 3-Sasakian orbifold $V$-bundle $\cals\ra{1.3}\calo$ over it 
\cite{BGM94}. 
Grove, Wilking and Ziller showed that in the case of Hitchin's metrics the bundle
$N^7_k\ra{1.3}M^4_k$ is actually a smooth 7-manifold and, 
moreover, computed $H^*(N^7_k,\bbz)$. When $k=2m-1$ the manifold $N^7_k$
turns out to be a rational homology 7-sphere with $H_3(N^7_k,\bbz)=\bbz_m$.

In dimension 5,
aside from the standard 5-sphere $S^5$ there is only one known case of an Einstein metric 
on a simply connected rational homology 5-sphere. It is the homogeneous space 
$SU(3)/SO(3)$ which is a simply connected non-spin manifold with 
$H_2(SU(3)/SO(3),\bbz)=\bbz_2.$ On the other hand the 5-manifolds $M^5$ considered in 
this paper are all spin, i.e. $w_2(M)=0.$ A well known theorem of Smale \cite{Sm62} says that 
for any simply connected 5-manifold with spin, the torsion group in $H_2$ is of the form 
$G\oplus G$ for some Abelian group $G.$ Until recently \cite{BG02} this was the only simply 
connected nontrivial rational homology sphere known to admit Riemannian metrics of positive 
Ricci curvature. In this note we go much further by proving the existence of Sasakian-Einstein 
metrics on infinitely many simply connected rational homology 5-spheres. Furthermore, the 
metrics typically depend on parameters, that is, there is non-trivial moduli.

We prove the following two theorems.

\begin{thm}\label{thm1}
There exists continuous parameter families of \Se metrics on infinitely many simply connected 
rational homology spheres in every odd dimension greater than 3. For some of these families 
the number of effective parameters grows exponentially with dimension.
\end{thm}

Dimension $5$ is treated separately since there is a classification of simply connected 
$5$-manifolds due to Smale \cite{Sm62} for spin manifolds and Barden \cite{Bar65} more 
generally, and our results are somewhat sharper. It is the spin manifold case that is relevant 
here, as any simply connected \Se manifold is necesarily spin. Our results
cover all orders of the torsion group $H_2(M_k,\bbz),$ that is not $4$ nor 
divisible by $6.$ We have

\begin{thm}\label{thm2}
For every integer $k>2$ that is either relatively prime to $3$ or $2$, there exist \Se 
metrics, depending on two real parameters, on a simply connected rational homology 
5-sphere $M_k^5$ with $w_2(M_k^5)=0,$ and $H_2(M_k^5,\bbz)$ having order $k^2$.  

\end{thm}

This result is probably not optimal. Furthermore, there are many more examples than those 
given in Section \ref{rathom5}.

\section{Branched Covers}

Let $f=f(z_1,\cdots,z_m)$ be a quasi-smooth weighted homogeneous polynomial of degree 
$d_f=w(f)$ in $m$ complex variables, and let $L_f$ denote its link. Let 
$\bfw_f=(w_1,\cdots,w_m)$ be the corresponding weight vector, and for a given weight vector 
$\bfw$ we let $|\bfw|=\sum_{i=1}^mw_i$ be its norm.
We consider branched covers constructed as the link $L_F$ of the polynomial
$$F=z_0^k+f(z_1,\cdots,z_m).$$
Then $L_F$ is a $k$-fold branched cover of $S^{2m+1}$ branched over the link $L_f.$ The 
degree of $L_F$ is $d_F=\lcm(k,d_f),$ and the weight vector is 
$\bfw_F=(\frac{d_f}{\gcd(k,d_f)},\frac{k}{\gcd(k,d_f)}\bfw_f).$ We shall always assume that 
$k\geq 2$, since the linear case $k=1$ is a hyperplane in a weighted projective space.
The following is Theorem 7.1 of \cite{BGN03b}:

\begin{thm}\label{torthm}
Let $f(z_1,\cdots,z_m)$ be a weighted 
homogeneous polynomial of degree $d$ and weights $\bfw=(w_1,\cdots,w_m)$ in 
$\bbc^m$ with an isolated singularity at the origin. Let $k\in \bbz^+$ and 
consider the link $L_g$ of the equation 
$$g=z_0^k+f(z_1,\cdots,z_m)=0.$$
Write the numbers ${d\over w_i}$ in irreducible form ${u_i\over v_i},$ and 
suppose that $\gcd(k,u_i)=1$ for each $i=1,\cdots,m.$ Then the link $L_g$ has weights 
$\frac{(d,k\bfw)}{\gcd(k,d)}$ and 
degree $\lcm(k,d).$ Furthermore, $L_g$  is 
a rational homology sphere with the order $|H_{m-1}(L_g,\bbz)|=k^{b_{m-2}(L_f)}$ where 
$b_{m-2}(L_f)$ is the $(m-2)$-nd Betti number of $L_f.$
\end{thm}

When considering rational homology spheres we can assume without loss of generality that 
$\gcd(k,d)=1.$ For if  $\gcd(k,d)>1$ and $\gcd(k,u_i)=1$ for each $i=1,\cdots,m,$ then any 
common factor of $k$ and $d$ must divide $w_i$ for each $i$, and so will be an overall 
common factor which gives an equivalent link. So hereafter we shall assume that 
$\gcd(k,d)=1.$

The formula for computing the $(m-2)$-nd Betti number can be obtained from Milnor and 
Orlik \cite{MO70}, viz.
\begin{equation}\label{Bettiequation}
b_{m-2}(L_f)= \sum (-1)^{m-s}\frac{u_{i_1}\cdots u_{i_s}}{v_{i_1}\cdots v_{i_s}{\rm 
lcm}(u_{i_1},\cdots,u_{i_s})},
\end{equation}
where the sum is taken over all the $2^{m}$ subsets $\{i_1,\cdots,i_s\}$ of $\{1,\cdots,m\}.$
This general form for $b_{m-2}$ is suitable for computer computations. 

\section{Sasakian-Einstein links and K\"ahler-Einstein orbifolds}

Let $X=X_{k,d,\bfw}$ denote the quotient orbifold by the natural $S^1$ action.
By now it is well understood that a \Se metric on the link $L_F$ is equivalent to a  
K\"ahler-Einstein metric of positive scalar curvature on  $X=X_{k,d,\bfw}.$  See \cite{BG01} 
for details.
It is easy to see that $X_{k,d,\bfw}$ is Fano if and only if 
\begin{equation}\label{Fano}
k(|\bfw_f|-d_f)+d_f>0.
\end{equation}
In particular, if $|\bfw_f|-d_f\geq 0,$ then $X_{k,d_f,\bfw}$ is Fano for all positive integers $k.$
There are three cases to consider. The case when the smaller link $L_f$ is itself Fano, i.e. 
$|\bfw_f|-d_f>0,$. This will be ruled out by Proposition \ref{notklt} below. The case when 
$L_f$ is Calabi-Yau, i.e. $|\bfw_f|-d_f= 0.$ This case is interesting since we have good 
solutions for infinitely many $k.$ Finally, there is the ``canonical'' case when $|\bfw_f|-d_f<0,$ 
for which the link $L_F$ is Fano for only finitely many $k.$ Nevertheless, there could be 
infinitely many good solutions in this case by taking arbitrarily large degrees $d_f.$

A sufficient condition for the existence of a K\"ahler-Einstein metric on the orbifold 
$X=X_{k,d_f,\bfw}$ is that for every effective $\bbq$-divisor $D$ that is numerically equivalent 
to $K^{-1}_X,$ the pair $(X,\frac{m-1+\gre}{m}D)$ is Kawamata log terminal or klt for some 
$\gre>0.$ For a precise definition see \cite{koll-mor}, and for a more complete discussion in 
the present context see \cite{dem-koll,BGN03,BGK03}. In particular, this implies the
condition 
\begin{equation}\label{klt}
k(|\bfw_f|-d_f)+d_f<\frac{m}{m-1}{\rm min}\{d_f,kw_i~|~i=1,\cdots m\}.
\end{equation}
This condition is necessary to satisfy the klt condition defined above, but it is far from 
sufficient, and we use it only to eliminate some cases. It can also be reiterated that while the 
klt condition on the pair $(X,\frac{m-1+\gre}{m}D)$ is sufficient to guarentee the existence of 
a K\"ahler-Einstein metric, it is far from necessary.

\begin{prop}\label{notklt}
If $\bfw_f-d>0$ then $X_{k,d,\bfw}$ does not satisfy the klt condition \ref{klt}.
\end{prop}

\begin{proof}
If $\bfw_f-d>0$ then the left hand side of \ref{klt} is at least $k+d.$ So for \ref{klt} to hold we 
must have
$$ (m-1)k<d, ~\text{and} \qquad (m-1)d<k.$$
This gives a contradiction.
\end{proof}

Next we look at the case $\bfw_f-d=0.$ Then the klt condition \ref{klt} becomes
\begin{equation}\label{klt2}
(m-1)d<mk~{\rm min}\{w_i\}.
\end{equation}
This is clearly satisfied for $k$ large enough, namely for $k>\frac{m-1}{m}d~{\rm 
max}\frac{1}{w_i}.$

We are particularly interested in the case of perturbations of Brieskorn-Pham (BP) 
singularities. 
In \cite{BGK03} J\'anos Koll\'ar and the authors gave sufficient conditions for the existence of 
K\"ahler-Einstein metrics on Fano orbifolds arising as perturbations of BP singularities. 
However, these conditions, although better than \ref{klt} above, are still not optimal.
We consider weighted homogeneous polynomials of the form
$$
P(z_0,\cdots,z_m)=\sum_{i=0}^mz_i^{a_i}+tp(z_0,\dots,z_m),
$$
where $t\in\bbc,$ and $w(p)=\lcm(a_0,\cdots,a_m)$ which is the degree of the polynomial. We 
impose a conditon on the zero set $Y(\bfa,p):=P^{-1}(0),$ namely 
the genericity
condition, which is always satisfied by $p\equiv 0$,
\begin{enumerate}
\item[(GC)] The intersections of $Y(\bfa,p)$ with any number of hyperplanes
$(z_i=0)$ \\
are all smooth outside the origin.
\end{enumerate}
Any polynomial or singularity satisfying condition GC is referred to as a (weighted 
homogeneous) {\it perturbation} of a Brieskorn-Pham polynomial or singularity.
Furthermore, 
we define
$$
C^j=\lcm(a_i:i\neq j),\quad
b_j=\gcd(a_j,C^j).
$$
The theorem proved in \cite{BGK03} is;

\begin{thm}\label{BP.KE.thm}
 The orbifold $Y(\bfa,p)/\c^*$ is Fano and 
has a \Ke metric if it satisfies condition (GC) and
$$
1<\sum_{i=0}^m\frac1{a_i}<
1+\frac{m}{m-1}\min_{i,j}\Bigl\{\frac1{a_i}, \frac1{b_ib_j}\Bigr\}.
$$
\end{thm}

\section{Rational Homology 5-Spheres}\label{rathom5}

We now specialize to the case of dimension 5, i.e. $m=3.$
Our 5-manifolds are constructed as $k$-fold branched covers of $S^5$ branched over 
certain Seifert manifolds that are in turn $S^1$ orbifold  V-bundles over a compact 
Riemann surface of genus $g.$ Our construction is similar to that in \cite{Sav79}. Let 
$f_3(z_1,z_2,z_3)$ be a 
weighted homogeneous polynomial of an isolated hypersurface singularity in $\bbc^3$ with 
weights $\bfw=(w_1,w_2,w_3)$ and degree $d.$ The link $L_\bfw$ defined by 
$L_\bfw=\{f_3=0\}\cap S^5$ is a Seifert fibration over an algebraic curve $C_\bfw$ in the 
weighted projective space $\bbp(\bfw).$ Let $g=g(\bfw)$ denote the genus of the curve 
$C_\bfw.$ In the case of rational homology 5-spheres, Theorem \ref{torthm} specializes to 
$|H_2(L_f,\bbz)|=k^{2g},$ while the Betti number formula \ref{Bettiequation} specializes to  
the genus formula of Orlik and Weigreich \cite{OW71}
$$g(C_\bfw)={1\over 2}\Bigl({d^2\over w_1w_2w_3}-d\sum_{i<j}{\gcd(w_i,w_j)\over 
w_iw_j}+\sum_i{\gcd(d,w_i)\over w_i}-1\Bigr).$$

The following lemma whose proof is clear from the genus formula will prove to be useful.
\begin{lem}\label{genuslem}
If $|\bfw|=d, \gcd(w_i,d)=w_i,$ and $\gcd(w_i,w_j)=1$ for all $i\neq j,$ then $g=1.$
\end{lem}

In view of Proposition \ref{notklt} we restrict ourselves to the case $\bfw_f-d\leq 0.$ In the 
case of Brieskorn 3-manifolds $L_f$ Milnor \cite{Mil75} refers to the cases 
$\bfw_f-d>0,\bfw_f-d=0,$ and $\bfw_f-d<0$ as spherical, Euclidean, and hyperbolic, 
respectively. In \cite{OW71} (with an erratum in \cite{OW77}) Orlik and Wagreich give a 
classification of all weighted homogeneous polynomials with only an isolated singularity at the 
origin up to $\bbc^*$-equivariant diffeomorphism. They gave a list of six classes, but there 
was an error and two more were reported in \cite{OW77} which completes their classification.
Here we analyze the possible $k$-fold branched covers when the Siefert manifold is 
Euclidean, i.e. $\bfw_f=d,$ which is straightforward. It turns out that the singularities arising 
here are precisely the unimodal parabolic singularities of \cite{AGV85}.  In any case it is easy 
to prove:

\begin{prop}\label{brEuclid}
Let $L_f$ be the link of a weighted homogeneous polynomial $f=f(z_1,z_2,z_3)$ of degree 
$d$ in three complex variables. Suppose further that $|\bfw|=d.$ Then $f$ is a weighted 
homogeneous polynomial with weight vector $\bfw,$ degree $d$, and number of monomials 
$n$ of degree $d$ given by one of the three cases:
\begin{center}\vbox{\[\begin{array}{|c|l|l|} \hline \bfw
& d &n\hl\hline 
(1,2,3)&6&7\hl
(1,1,2)&4&9\hl
(1,1,1)&3&10\hl
\end{array}\]}
\end{center} 
\end{prop}

Let $L_F$ be a $k$-fold branched cover of a 3-dimensional Euclidean link, then by 
Proposition \ref{brEuclid} $L_F$ is the link of an isolated hypersurface singularity of the form
$$z_0^k+f_d(z_1,z_2,z_3)=0$$
with weights $(d,k\bfw_f)$ where $d$ and $\bfw_f$ are given in the table above, and $f_d$ is 
an arbitrary polynomial polynomial of the given weights subject to the condition that the 
singularitiy is isolated. There is a dense open subset of the parameter space where this is a 
perturbation of a Brieskorn-Pham singularity. It is easy to see that the bounds given in 
Theorem \ref{BP.KE.thm} are satisfied if we choose $k\geq 3$ in the $d=4$ and $d=3$ cases, 
and $k\geq 5$ in the $d=6$ cases. Hence, in these cases we get \Se metrics on the links 
$L_F.$ Furthermore, in the $d=3$ and $d=6$ cases there is a two real parameter family of 
solutions. This proves Theorem \ref{thm2}. \qed

Proposition \ref{brEuclid} gives a classification of all the possible Euclidean cases. We can also 
consider hyperbolic cases. Since the weights are unrestricted there are many, so we only 
mention one which is a special case of Example \ref{brcan} below. We consider links 
associated to perturbations of the BP singularity
$$z_0^k+z_1^l+z_2^l+z_3^l=0.$$
The Betti number formula \ref{Bettihyper} below simplifies to
$b_1=(l-2)(l-1).$ Thus, by Theorem \ref{torthm} we get rational homology 5-spheres 
$M^5_{k,l}$ with $H_2(M^5_{k,l},\bbz)$ of order $k^{(l-2)(l-1)}.$ From equation \ref{lbound} 
below we see that $l$ must be either $4$ or $5.$ The Fano condition and klt inequalities of 
Theorem \ref{BP.KE.thm} are satisfied only when $k=3$ and $l=4.$
By Equation \ref{modhyper} the number of effective real parameters is $12.$ Thus, there is a 
rational homology 5-sphere $M^5_{3,4}$ with $|H_2(M^5_{3,4},\bbz)|=3^6$ that admits a \Se 
metric depending on $12$ real parameters.

\medskip

\section{Examples in Dimension $2m-1$}

\begin{exmp}\label{Fermat}[Branched Covers of Calabi-Yau hypersurfaces]
First we consider $k$-fold branched covers of Fermat-Calabi-Yau hypersurfaces
\begin{equation}\label{brFermat}
z_0^k+z_1^m+\cdots z_m^m=0
\end{equation}
with $\gcd(k,m)=1.$  For $m\geq 3$ the link $M^{2m-1}_k$ is a simply connected nontrivial 
rational homology sphere. The $(m-2)$-nd Betti number of the Calabi-Yau link is
$$b_{m-2}= (-1)^m\Bigl(1+\frac{(1-m)^m-1}{m}\Bigr).$$
So by Theorem \ref{torthm} $|H_{m-1}(M^{2m-1}_k,\bbz)|=k^{b_{m-2}}.$ For example, in 
dimension $7$ ($m=4$) we get $|H_2(M^7_k,\bbz)|=k^{21}$ whereas, in dimension $9$ we 
have $|H_2(M^9_k,\bbz)|=k^{204}.$

It is easy to see that the klt condition in Theorem \ref{klt} is satisfied if $k>m(m-1)$ in which 
case the rational homology spheres $M^{2m-1}_k$ admits a family of \Se metrics. The 
number of effective complex parameters $\mu$ is determined by (cf. \cite{BGK03})
$$\mu=h^0(\bbp(\bfw),\calo(d))-\sum_ih^0(\bbp(\bfw),\calo(w_i))$$
where $\bfw=(m,k,\cdots,k)$ and $d=mk.$ Thus, we find 
$$\mu={2m-1 \choose m}-m^2.$$ 
By Sterling's formula one sees that $\mu$ grows exponentially with $m.$
In dimension $7$ we get $38$ real parameters, while in dimension $9$ we have $202$ 
effective real parameters. This example is enough to prove Theorem \ref{thm1}. \qed

In the orbifold category there are many Calabi-Yau hypersurfaces in weighted projective 
spaces $\bbp(\bfw).$ For example, for general $m$, one can consider $\bfw=(1,\cdots, 
1,\break m-1)$ of degree $d=2(m-1),$ or $\bfw=(1,\cdots,1,m-2,m-2)$ of degree 
$d=3(m-2),$ etc.
For rational homology spheres of dimension 7, we obtain examples from branched covers of 
Reid's list (cf. \cite{Fle00}) of 95 log K3 surfaces, and in dimension 9 from branched covers of 
the over 6000 Calabi-Yau orbifolds in complex dimension 3 \cite{CLS90}.

\end{exmp}

\begin{exmp}[Branched Covers in the Canonical Case]\label{brcan}
We consider links $L_F$ of branched covers of a canonical Fermat hypersurface of the form
$$F=z_0^k+z_1^l+\cdots +z_m^l=0$$
with $l\geq m+1$ and $\gcd(k,l)=1.$ Here the quotient $X$ will be Fano if and only if 
$$k <\frac{l}{l-m}.$$
Also for any branched cover we are only interested in $k\geq 2.$ Combining this with the Fano 
condition gives an upper bound on $l$ for fixed $m,$ namely $l< 2m.$ Thus, we get the range 
for $l$
\begin{equation}\label{lbound}
m+1\leq l\leq 2m-1.
\end{equation}
The link $L_F=M^{2m-1}_k$ is a rational homology sphere with 
$|H_{m-1}(M^{2m-1}_k,\bbz)|=k^{b_{m-2}}$ where
\begin{equation}\label{Bettihyper}
b_{m-2}= (-1)^m\Bigl(1+\frac{(1-l)^m-1}{l}\Bigr).
\end{equation}

We now look at the klt condition, that is, the right hand inequality of Theorem \ref{BP.KE.thm}. 
We see that $L_F$ will have \Se metrics if we choose $k$ to satisfy
\begin{equation}\label{hyperklt}
\frac{(m-1)l^2}{(m-1)l(l-m)+m}< k<\frac{l}{l-m}.
\end{equation}

As above we can compute the number of effective complex parameters $\mu$ to be
\begin{equation}\label{modhyper}
\mu={m+l-1 \choose l}-m^2.
\end{equation}
As a special case we consider $l=m+1$. In this case the only solution to the inequality 
\ref{hyperklt} is $k=m.$ This gives \Se metrics on rational homology spheres depending on 
$2({2m \choose m+1}-m^2)$ effective real parameters.

Other solutions can be worked out, for example, the singularity
$$z_0^k+z_1^{2m}+\cdots +z_{m-1}^{2m}+z_m^2=0.$$
This satisfies the inequalities of Theorem \ref{BP.KE.thm} if we choose $k=2m-1.$

\end{exmp}

\begin{ack}  
We would like to thank J\'anos Koll\'ar for helpful discussions concerning this paper.
The authors were partially supported by the NSF under grant number
DMS-0203219. 
\end{ack}

\def\cprime{$'$} \def\cprime{$'$}
\providecommand{\bysame}{\leavevmode\hbox to3em{\hrulefill}\thinspace}
\providecommand{\MR}{\relax\ifhmode\unskip\space\fi MR }
\providecommand{\MRhref}[2]{%
  \href{http://www.ams.org/mathscinet-getitem?mr=#1}{#2}
}
\providecommand{\href}[2]{#2}


\end{document}